\documentclass[12pt]{article}
\usepackage[german,english]{babel}
\usepackage{amsfonts}
\usepackage{amssymb}
\usepackage{amsbsy}
\usepackage{amscd}
\usepackage{graphicx}
\usepackage{amsmath}
\usepackage{fabfeynmp,subfig}
\input epsf.sty

\newcommand{\BEQ}{\begin{equation}}     % Gleichungen Anfang ..
\newcommand{\BEA}{\begin{eqnarray}}
\newcommand{\EEQ}{\end{equation}}       % .. und Ende
\newcommand{\EEA}{\end{eqnarray}}
\newcommand{\g}{{\mathfrak{g}}}

\newcommand{\h}{{\mathfrak{h}}}

\newcommand{\R}{\mathbb{R}}

\newcommand{\C}{\mathbb{C}}
\newcommand{\Z}{\mathbb{Z}}

\def\abs#1{\left| #1 \right|}

\def\schwx#1,#2{{#1'''(#2) \over #1'(#2)} - {3 \over 2}{#1''(#2)^2 \over #1'(#2)^2}}

\def\norm#1{\left\| #1 \right\|}%%%%%%%%                                   % Nummerierung pro section
\begin{document}
%\begin{fmffile}{ARTHISTMpost}
\begin{titlepage}
\begin{center}
{\Large {\bf Sur les origines du cocycle de Virasoro}}
\end{center}
\vskip5cm

\centerline{ {\bf Claude Roger}} 
\vskip 0.5 cm
\begin{center}
  {$^a$Institut Camille Jordan ,\footnote{Laboratoire associ\'e au CNRS UMR 5208} 
Ecole Centrale de Lyon, INSA de Lyon, Universit\'e de Lyon, Universit\'e
  Lyon I,} 
  \end{center}
  \centerline
{43 boulevard du 11 novembre 1918, 
F-69622 Villeurbanne Cedex, France}
\vskip5cm

\begin{abstract}
\noindent
Cet article retrace un bref historique des origines du cocycle de Virasoro, en alg\`ebre et en  th\'eorie quantique des champs.
\vskip1cm
\noindent
This article gives a short sketch of the origins of Virasoro cocycle, both in algebra and quantum field theory.
\vskip1cm
\noindent
{\it NB: cet article est une version longue de l'appendice historique de \cite{GuiRog}, auquel le lecteur pourra se r\'ef\'erer pour plus de d\'etails sur les notions introduites dans la suite.}
\end{abstract}

\end{titlepage}
\begin{fmffile}{ARTHISTMpost}
\large
Cet article a pour but de retracer l'histoire du cocycle de Virasoro et de certains analogues qui lui
sont indissociablement li\'es; nous n'avons pas cherch\' e ici \`a faire r\'eellement oeuvre d'historien, ce qui serait pr\'ematur\'e, nous avons simplement tent\'e  de d\'ecrire la g\'en\'ealogie de cet objet physico-math\'ematique,  de retrouver ses
racines tant math\'ematiques que physiques, sans toutefois faire le bilan de tous les endroits o\`u ce cocycle appara\^it; nous aurons ainsi \`a
parcourir des domaines aussi vari\'es que fascinants, \`a premi\`ere
vue tr\`es \'eloign\'es les uns des autres, de l'alg\`ebre homologique
\`a la th\'eorie des champs, et nous  y verrons certaines constructions alg\'ebriques ou g\'eom\'etriques revenir p\'eriodiquement dans des domaines a priori diff\'erents.

\smallskip

%%%%%%%%%%%
Mais rappelons tout d'abord de quoi il s'agit: du cercle, de ses
diff\'eomorphismes, et de ses champs de vecteurs tangents. Tout champ
tangent au cercle s'identifie \`a une fonction sur le cercle gr\^ace au
choix d'un param\'etrage, et le crochet de Lie des champs de vecteurs peut alors
s'\'ecrire :
$$
[f(t){d\over dt}, g(t){d\over dt}]=(f(t)g'(t)-g(t)f'(t)){d\over dt}
$$
Cette alg\`ebre de Lie est not\'ee ${\rm Vect}(S^1)$. L'utilisation de la base de Fourier $e_n=i e^{int}{d\over
dt}$,  dans laquelle le crochet s'\'ecrit :
$$
[e_n, e_m]=(n-m)\ e_{n+m}.
$$ permet d'en donner une description alg\'ebrique qui se g\'en\'eralise naturellement \`a coefficients dans un corps quelconque de caract\'eristique z\' ero.
Cette alg\`ebre de Lie  prend le nom d'{\it alg\`ebre de Virasoro}
quand elle se voit rajouter un myst\'erieux terme central :
\begin{equation}
\displaystyle 
[f(t){d\over dt}, g(t){d\over dt}]
=
(f(t)g'(t)-g(t)f'(t)){d\over dt} + c \int_{S^1}f(t)g'''(t)dt
\tag{1}
\end{equation}
soit, en termes des g\'en\'erateurs, la  formule devenue c\'el\`ebre:
\begin{equation}
\displaystyle
[e_n, e_m]=(m-n)e_{n+m} + {n^3\over 12} \delta_{0,n+m} c. 
\tag{2}
\end{equation}

\smallskip

A chaque alg\`ebre de Lie son groupe? Malheureusement, ce n'est plus toujours vrai en dimension infinie, mais ici se produit un petit miracle, du moins si nous restons en coefficients r\'eels:
le groupe ${\rm Diff}(S^1)$ des diff\'eomorphismes du cercle peut \^etre vu
comme
 un groupe de Lie  de dimension infinie dont l'alg\`ebre de Lie s'identifie
\`a ${\rm Vect}(S^1)$; de plus lui aussi peut \^etre agr\'ement\'e d'un terme
central !
En voici la formule:

on  peut  tout d'abord associer \`a tout $\displaystyle f \in {\rm Diff}(S^1)$ une fonction
sur ${S^1}$ not\'ee $\mu_f$ d\'efinie par $\displaystyle
f^*(\omega)=\mu_f \omega$. Le cocycle central, dit de Bott-Thurston, se d\'efinit ensuite par la
formule
\begin{equation}
\displaystyle
BT(f,g)=\int_{S^1} \log \mu_{f \o g}d\log \mu_g .
\tag{3}
\end{equation}
%%%%%%%%%
On peut distinguer sommairement deux racines principales de ce curieux objet math\' ematique: l'une est strictement alg\'ebrique avec les {\it multiplicateurs de Schur}\cite{Schur} qui
marqu\`erent les d\'ebuts de l'alg\`ebre homologique, l'autre physique
avec les c\'el\`ebres {\it termes de Schwinger}, premi\`ere apparition d'un ph\'enom\`ene d'
"anomalie" en th\'eorie quantique des champs. Le point commun entre ces
deux origines est l'ubiquit\'e d'un  autre cocycle, universel celui-l\`a, qui se pr\'esente sous la forme
$c(A,B)={\rm Tr} (A \delta B)$ (pour les d\'etails, cf.infra) comme celui qui gouverne l'extension centrale des matrices de Jacobi \cite{KaRa} ou encore les groupes lin\'eaires
restreints \cite{PressleySegal}. Nous le rencontrerons tout au long de cet article; il est \`a la source, non seulement du
cocycle de Virasoro, mais de bien d'autres analogues comme ceux des groupes de jauge ou ceux des
alg\`ebres de Kac-Moody; rappelons pourquoi ces derni\`eres sont indissociablement li\'es \`a l'alg\`ebre de Virasoro. Les alg\`ebres de Kac-Moody  sont pr\'ecis\'ement obtenues par extensions \`a partir des  alg\`ebres de lacets , et celles ci sont des alg\`ebres de Lie de fonctions sur $S^1$ \`a valeurs dans des alg\`ebres de Lie semi-simples de dimension finie. L'action du groupe $ {\rm Diff}(S^1)$ par reparam\'etrisation est alors naturelle...Pour plus de d\'etails, nous renvoyons \`a \cite{GuiRog}.

\smallskip
%%%%%%%%%%%%%%%%%
L'aspect purement alg\'ebrique est celui de la th\'eorie des extensions
centrales des groupes et alg\`ebres de Lie, dont on peut trouver par exemple les d\'etails dans le trait\'e de Hilton et Stammbach \cite{HiltonStammbach}; les deux types d'extensions se traitent de fa\c con tout \`a fait parall\`eles, les technicit\'es cohomologiques \'etant variables suivant les cas. Il est temps de proc\'eder  maintenent \`a quelques rappels.
%%%%%%%%%%%%%%%%%%
Soient ${\g}$ et ${\h}$ deux alg\`ebres de Lie, 
une {\it extension de ${\g}$ par ${\h}$} est une alg\`ebre de Lie
$\Hat{\g}$ telle qu'il existe une suite exacte d'alg\`ebres de Lie
\begin{equation}
\displaystyle
\begin{CD}
0@>>>\h @>{i}>>\Hat{\g}@>{\pi}>> \g@>>> 0\end{CD}
\tag{4}
\end{equation}
Une extension de ${\g}$ par $\h$ est {\it centrale} si 
$i(\h)$ est contenu dans le centre de $\Hat{\g}$. Les extensions centrales de groupes, discrets ou de Lie, se d\' efinissent de fa\c con analogue. Mais pr\'ecisons maintenant une des raisons de leur importance, \`a savoir les repr\'esentations projectives:

Une {\it repr\'esentation projective du groupe $G$} est une
application $\rho:G\to GL(n,\mathbb{K})$ telle qu'il existe $$
\begin{CD}
G\times G @>{c}>>\mathbb{K}^*\end{CD}$$ 
telle que l'on ait
 pour tous les $x,y$ dans $G$,$$\rho(x)\rho(y)=c(x,y)\rho(xy).$$ En d'autres termes on a un homomorphisme \`a valeurs dans le groupe projectif et non le groupe lin\'eaire; une repr\'esentation projective peut se lin\'eariser en une repr\'esentation lin\'eaire d'une extension centrale ad-hoc.

%%%%%%%%%%%%%%%%%
%%%%%%%%%%%%%%%
\bigskip

Au niveau des groupes discrets, rappelons le r\'esultat fondamental suivant : si $G$ est un groupe
parfait, c'est-\`a-dire si $H_1(G)=G/[G,G]=0$, alors il existe une
{\it extension centrale universelle}. 

\begin{equation}
\displaystyle
\begin{CD}
1 @>>>H_2 (G)@>>>  \Hat{G}@>>>G @>>> 1\end{CD}
\tag{5}
\end{equation}
Cette notion d'universalit\' e est \`a comprendre dans le sens
cat\'egorique usuel d'objet initial dans la cat\'egorie des extensions centrales, chacune \' etant quotient de l'extension universelle. 
Plus g\'en\'eralement toutes les extensions centrales $(1) \to C \to E
\to G \to (1)$ sont classifi\'ees par une classe de cohomologie $c \in
H^2 (G,C)=Hom (H_2 (G), C)$. On peut \'etablir en outre une bijection
entre les classes d'\'equivalence d'extensions centrales de $G$ et les
sous groupes du groupe ab\' elien $H_2(G)$ en faisant correspondre \`a chaque extension
le sous-groupe $C \subset H_2 (G)$. Ce th\'eor\`eme est d\^u sous
cette forme tr\`es g\'en\'erale \`a Michel Kervaire (1970) \cite{Ke}, mais
une premi\`ere version en est apparue d\`es 1904 pour le cas o\`u le
groupe $G$ est fini, le groupe $H_2 (G)$ correspondant alors 
aux {\it multiplicateurs de Schur}; des groupes finis \`a ${\rm Diff}(S^1)$
le chemin a \'et\'e long et panoramique, mais le cadre alg\'ebrique est rest\'e presque identique.

Une \'etape interm\'ediaire
essentielle fut accomplie par Heinz Hopf avec son article de 1942 \cite{Hop}
contenant la c\'el\`ebre pr\'esentation du $H_2$ par g\'en\'erateurs
et relations: $H_2 (G)= (R \cap [F,F])/[R ,F]$. Ici $R$ repr\'esente le noyau d'une surjection
$F\to G$ o\`u $F$ est un groupe libre(voir \cite{HiltonStammbach}). Dans ce m\^eme travail, on trouve
  une remarque  qui devait s'av\'erer aussi profonde que fructueuse pour l'avenir
: l'analogie entre les extensions centrales de groupes et les
rev\^etements des espaces topologiques, $H_2 (G)$ correspondant alors
au groupe fondamental $\pi_1 (X)$. ... Curieusement, la th\'eorie des extensions centrales pour
les alg\`ebres de Lie dut attendre 1982 (Loday et Kassel \cite{KaLo}) pour
\^etre publi\'ee explicitement, m\^eme si elle \'etait consid\'er\'ee
comme r\'esultat ``folklorique''. Le probl\`eme de l'extension centrale universelle pour une alg\`ebre de Lie se traite
comme pour les groupes : si $H_1 (\g)=0$, alors $\g$ admet une
extension centrale universelle 

\smallskip
\begin{equation}
\displaystyle
\begin{CD}
0 @>>>H_2 (\g)@>>>  \Hat{\g}@>>>\g @>>> 0\end{CD}\tag{6}
\end{equation}

La cohomologie des alg\`ebres de
Lie fut d\'efinie par Chevalley et Eilenberg \cite{ChEi} puis d\'evelopp\'ee
par J. L. Koszul \cite{Kos} dans sa th\`ese d'Etat vers 1950; on peut toutefois identifier sa naissance, sous une forme crypt\'ee, dans une note d'Elie Cartan en 1928\cite{Cartan}, dans lequel il consid\`ere le sous complexe des formes diff\'erentielles invariantes sur un groupe de Lie compact. Bien plus tard, cette note devait inspirer D.Sullivan dans ses travaux sur l'homotopie rationnelle, mais ceci est une autre histoire...

 Ce ne fut que
vers les ann\'ees 60 que la cohomologie des alg\`ebres de Lie connut son plein d\'eveloppement  avec
ses applications en th\'eorie des repr\'esentations, topologie
diff\'erentielle, et d\'ej\`a physique th\'eorique comme nous allons
le voir bient\^ot; le magistral trait\'e de H.Cartan et S.Eilenberg\cite{CarEil}  y a certainement jou\'e un r\^ole de catalyseur. On se limitera ici \`a rappeler le r\^ole des cohomologies continues
pour les repr\'esentations des groupes ainsi que celui des cohomologies de certaines alg\`ebres de Lie pour l'\'etude des classifiants et classes caract\'eristiques.

C'est \'egalement vers cette \'epoque que commenc\`erent les travaux
de Gelfand et Fuks sur la cohomologie des alg\`ebres de Lie de champs
formels et des champs de vecteurs. C'\'etait une p\'eriode brillante des math\'ematiques moscovites, que certains aiment encore \`a d\'esigner sous le nom d' "\'ecole russe" (mais pourquoi ignorer ceux de P\'etersbourg  sous la houlette des Faddeev p\`ere et fils ?); le fameux s\'eminaire Gelfand vit s'\'elaborer  parmi beaucoup d'autres choses la m\'ethode des orbites en th\'eorie des repr\'esentations, les cohomologies dites de Gelfand-Fuks, la superg\'eom\'etrie et maintes applications alg\'ebriques, g\'eom\'etriques, et physiques; cette ouverture des math\'ematiques suppos\'ees pures vers la physique (et r\'eciproquement) commen\c cait \`a apporter un d\'ementi au divorce d\'eplor\'e par Dyson dans sa conf\'erence \`a l'AMS en 1972; surtout, le contraste avec la situation fran\c caise de la m\^eme \'epoque \'etait spectaculaire !

Avant d'aller plus loin, pr\'ecisons les raisons de s'int\'eresser \`a ces alg\`ebres de Lie des champs de vecteurs sur les vari\'et\'es: ce sont les alg\`ebres de Lie correspondant aux pseudogroupes de Lie (et non groupes, qui n'existent pas vraiment dans ce contexte) de diff\'eomorphismes des ouverts des $\mathbb{R}^n$ 
ou des vari\'et\'es, leur classification aboutit \`a celle des diff\'erents types de "structures" sur les vari\'et\'es, la connaissance de leurs invariants(dont la cohomologie) est donc essentielle. Le calcul de $H^*({\rm
Vect}(S^1))$ contenant la premi\`ere apparition du $2$-cocycle central
fut publi\'e en 1968 \cite{FuGe1}, le cas des vari\'et\'es de dimension sup\'erieure fut \'elucid\'e peu apr\`es, et l'alg\`ebre de Lie ${\rm Vect}(S^1))$ (que les math\'ematiciens n'appelaient pas encore de Virasoro) apparut d\`es lors comme tr\`es particuli\`ere, car elle \'etait{\it la seule \`a admettre une extension centrale non triviale}, ses grandes soeurs de dimension $n>1$ v\' erifiant $H^2=0$.

Un autre r\'esultat va mettre en relief l'importance de
l'alg\`ebre de Virasoro parmi les alg\`ebres de Lie simples actuellement connues. Afin d' \'etendre \`a la dimension infinie la classification des
alg\`ebres simples
de dimension finie, on peut s'int\'eresser aux alg\`ebres
$\Z$-gradu\'ees, soit $ \g = \underset{n \in \Z}{\bigoplus} \g_{n}$
telles que le crochet v\'erifie $[\g_{n}, \g_{m}] \subset \g_{n+m}$, chacun des sous-espaces
$\g_{n}$ \'etant de
dimension finie
$p(n)$; on suppose que $\g$ est simple et que $p(n)$ est \`a croissance au
plus polyn\^omiale.
Parmi ces alg\`ebres de Lie, on connait bien :
\begin{enumerate}

\item {\it Les alg\`ebres de Lie semi-simples de dimension finie}, pour
leur graduation donn\'ee
par un syst\`eme de racines. On a alors $\g_{n} = \{0\}$ sauf pour un nombre
fini de valeurs de $n$.

\item {\it Les alg\`ebres de lacets} construites sur les alg\`ebres
semi-simples de dimension finie et leurs extensions centrales dites
\textit{alg\`ebres affines ou alg\`ebres de Kac-Moody}, d\'ej\`a mentionn\'ees dans l'introduction comme indissociables
de l'alg\`ebre de Virasoro.

\item {\it Les alg\`ebres de Lie dites de Cartan}: ce sont les
alg\`ebres de champs de vecteurs
formels ${\rm Vect}(n)$ et leurs sous-alg\`ebres associ\'ees aux alg\`ebres
de Lie de champs de
vecteurs simples,
primitives et transitives; parmi elles on trouve
 les alg\`ebres
 de Lie des champs de vecteurs formels, et leurs sous-alg\`ebres des champs
symplectiques,
unimodulaires, et de contact. Ce sont les analogues formels des alg\`ebres de champs de vecteurs tangents aux vari\'et\'es.
\item \textit{L'alg\`ebre de Virasoro} (dans sa version alg\'ebrique avec les g\'en\'erateurs de Fourier).
\end{enumerate}

Une conjecture de V.G. Kac affirmant que ces alg\`ebres 
\'etaient les seules
 \`a v\'erifier
les conditions de simplicit\'e, $\Z$-graduation, et de croissance au
plus polyn\^omiale de $p(n)$ a \'et\'e transform\'ee en th\'eor\`eme
gr\^ace \`a un spectaculaire r\'esultat d'Olivier Mathieu \cite{Mat1}. 

Il faut ici pr\'eciser un point de terminologie souvent aga\c cant: certains auteurs, principalement physiciens, ont tendance \`a appeler l'alg\`ebre de Virasoro sans terme central, surtout lorsqu'elle est donn\'ee sous sa forme alg\'ebrique par $
[e_n, e_m]=(n-m)e_{n+m},
$ l ' "alg\`ebre de Witt". Si l'on se reporte aux travaux de l'\'eminent alg\'ebriste et arithm\'eticien Ernst Witt
et ses collaborateurs\cite{Witt,Chang}, on s'aper\c coit vite qu'ils ne concernent que les cas o\`u la caract\'eristique est non nulle, et il vaut mieux r\'eserver le terme d'"alg\`ebre de Witt" \`a ce cas, qui pose d'ailleurs des probl\`emes cohomologiques nettement plus difficiles\cite{Leites}.

  Un sujet voisin des cohomologies dites de Gelfand -Fuks, et qui s'est bien d\'evelopp\'e au cours des ann\'ees 80, a \'et\'e l'\'etude 
de la cohomologie de l'alg\`ebre de Lie $\g_{\mathcal A}$ obtenue
par produit tensoriel de la $\mathbb{K}$-alg\`ebre de Lie $\g$ et de la
$\mathbb{K}$-alg\`ebre associative et commutative $\mathcal A$, le crochet \'etant
d\'efini par $[X \otimes a, Y \otimes b]=[X, Y]\otimes ab$. Les
alg\`ebres de courants, apparues en physique avec les th\'eories de
jauge non ab\'eliennes (champs de Yang-Mills)\cite{Zee}
rel\`event de cette
cat\'egorie avec $\mathbb{K}=\mathbb{R}$, $\mathcal A$ une alg\`ebre de fonctions
$C^\infty$ sur un espace-temps, et $\g$ une alg\`ebre de Lie
semi-simple, repr\'esentant les sym\'etries internes du probl\`eme, le
plus souvent $\g=\mathfrak{su}(n)$. Si on consid\`ere maintenant le
cas tr\`es particulier o\`u $\mathcal A =C^\infty (S^1)$, l'alg\`ebre
$\g_{\mathcal A}$ n'est autre que l'alg\`ebre de Lie des lacets sur $\g$
; l'extension centrale universelle de cette alg\`ebre de Lie a \'et\'e
mise en \'evidence ind\'ependamment par V. Kac \cite{Kac3} et R. Moody 
\cite{Mo}
dans leur imposante classification des alg\`ebres de Lie
filtr\'ees. Ces extensions centrales devaient conna\^\i tre une
brillante destin\'ee, sous le nom d'alg\`ebres de Kac-Moody;
remarquons cependant que dans son c\'el\`ebre livre, Victor Kac
indique au Chap 7...{\it``the formula for the cocycle has been known
to physicists for such a long time that it is difficult to trace the
original source...''}. Rappelons donc cette illustre formule :

Pour $f$ et $g$ dans $\g \otimes C^\infty (S^1)=C^\infty (S^1,
\g)$, on a \begin{equation}
\displaystyle c (f,g)=\int_{S^1} \kappa (f dg)
\tag{7}
\end{equation}
o\`u
$\kappa : \g \times \g \to \mathbb{K}$ d\'esigne la forme de Killing.

 Si l'on
pr\'ef\`ere une formulation complexe, on a pour $f$ et $g$ dans $\g
\otimes \C[z, z^{-1}]$, $\displaystyle c (f,g)={\rm Res}_{z=0} \kappa
(fdg)$. 

Sous ses deux formes, ce cocycle est obtenu \`a partir des
deux m\^emes ingr\'edients: d'une part une forme bilin\'eaire invariante, ici
la forme de Killing accompagn\'ee de l'int\'egrale, ou du r\'esidu ce
qui revient au m\^eme, d'autre part une d\'erivation ext\'erieure de
l'alg\`ebre de Lie $\g_{\mathcal A}$, ici l'application $f \mapsto f'$. En
toute rigueur, la preuve de l'universalit\'e de cette extension repose
sur le calcul de la cohomologie des alg\`ebres de lacets (voir par
exemple Pressley et Segal \cite{PressleySegal} ou Lepowsky \cite{Lep}).

Dans certains cas cette construction peut se
g\'en\'eraliser : si $\g$ est une alg\`ebre de Lie munie d'une forme
invariante $\kappa : \g \times \g \to\mathbb{K}$, on peut associer \`a une
d\'erivation $\delta : \g \to \g$, un $2$-cocycle $c_\delta : \g
\times \g \to \mathbb{K}$ sous certaines conditions d'antisym\'etrie; 
pour ce qui est de la cohomologie de $\g_{\mathcal
A}$, le pas d\'ecisif fut franchi par S. Bloch dans un article
m\'emorable \cite{Bl}. Si $\g=\mathfrak{sl}(n)$, on a $H_2 (\g_{\mathcal
A})=\Omega^1 \mathcal A/d \mathcal A$ o\`u $\Omega^1 \mathcal A$ d\'esigne
l'espace des diff\'erentielles de K\"ahler sur $\mathcal A$ et $d\mathcal A$
le sous-espace des diff\'erentielles exactes. Ce r\'esultat s'obtient
par la construction explicite d'une extension centrale 
%\begin{equation*}
%\label{A3}
%\xymatrix{0 \ar[r]& \Omega^1 \mathfrak A /d \mathfrak A \ar[r]& \wh{\g}_{\mathfrak
%A} \ar[r]& \g_{\mathfrak A} \ar[r]& 0}
%\tag{3}
%\end{equation}
\begin{equation}
\displaystyle
\begin{CD}
0 @>>> \Omega^1 \mathcal A/d \mathcal A @>>> \hat{ \g}_{\mathcal A}@>>>\g_{\mathcal A}@>>>0\end{CD}
\tag{8}
\end{equation}

dont Spencer Bloch a montr\'e l'universalit\'e, au moins pour 
$\g=\frak{sl}(n)$ avec $n$ suffisamment grand. Remarquons que si $\g_{\mathcal A}$ est l'alg\`ebre des courants sur une vari\'et\'e $V$, alors $\Omega^1
\mathcal A /d \mathcal A=\Omega^1 (V)/d \Omega^0 (V)$ (le calcul du $H_2$ a
\'et\'e r\'ealis\'e dans ce cas par B. Feigin\cite{Feigin},
ind\'ependamment du r\'esultat de S. Bloch). En particulier si
$V=S^1$, $\Omega^1 \mathcal A /d \mathcal A$ est alors isomorphe \`a $\R$
via l'int\'egrale et on retrouve l'alg\`ebre de Kac-Moody associ\'ee
\`a $\g$. D'ailleurs, dans l'introduction de \cite{Bl}, S. Bloch indique
dans son remerciement \`a Deligne que ...{\it ``In communicating to me his
construction, Deligne remarked (cryptically, as his wont) that he had
found it while thinking about Kac-Moody Lie algebras...''}. Une \'etape
importante de cette g\'en\'eralisation du cocycle de Kac-Moody a
\'et\'e la construction g\'en\'erale des traces et des r\'esidus pour
les diff\'erentielles sur des courbes quelconques, due \`a J. Tate
\cite{Tat}. Ce dernier r\'esultat peut \^etre consid\'er\'e comme une version pr\'eliminaire
des r\'esidus de dimension sup\'erieure, dus \`a Grothendieck et Hartshorne.

La
forme explicite du cocycle de l'extension pr\'ec\'edente est la suivante 
\begin{equation}
\displaystyle
c(X \otimes a, Y \otimes b)=\kappa (X,Y)[adb] \tag{9}
\end{equation}
o\`u $[adb]$ d\'esigne la classe de la diff\'erentielle $adb$ modulo
les diff\'erentielles exactes. Si on a une application $I : \Omega^1
\mathcal A /d \mathcal A  \to \mathbb {K} $, comme par exemple l'int\'egration sur un
cycle pour $\mathcal A=C^\infty (V)$, on en d\'eduit une extension
centrale \`a noyau scalaire par la formule $c_I (X \otimes a, Y
\otimes b)=\kappa (X,Y)I(adb)$.

Ce r\'esultat inspira les travaux de Loday et Quillen qui purent
calculer l'homologie de $\g_{\mathcal A}$ quand 
$\g=\mathfrak{gl}(\infty)$ $= \displaystyle \underset{\longrightarrow}{\rm lim}
\, \, \mathfrak{gl}(n)$. Cette homologie est munie naturellement d'une
structure d'alg\`ebre de Hopf, dont ils ont d\'etermin\'e l'espace des
\'el\'ements primitifs ; on obtient 
$${\rm Prim} (H_* (\g_{\mathcal A}, \mathbb{K}))=HC_* (\mathcal A)$$
avec un d\'ecalage d'indices (Thm. de Loday et Quillen\cite{Loday1}). La notation $HC_*$
d\'esigne ici l'homologie cyclique, et ce r\'esultat g\'en\'eralise
bien celui de Bloch car $HC_1 (\mathcal A)=\Omega^1 \mathcal A /d\mathcal A$ si
$\mathcal A$ est une $\mathbb{K}$-alg\`ebre commutative lisse; Feigin et Tsygan
 avaient construit cette homologie cyclique ind\'ependamment de A Connes, sous le nom de
``$K$-th\'eorie additive''\cite{FeTsy}, mais en quoi est-elle bien une version lin\'eaire ou additive de la $K$-th\'eorie?
Au terme  $H_2(\g_{\mathcal A})=HC_1 (\mathcal A)$
correspond $H_2 (E (\mathcal A))=K_2(\mathcal A)$ avec l'extension centrale
de groupes 
%\begin{equation*}
%\xymatrix{(1) \ar[r]& K_2(\mathfrak A) \ar[r]& {\rm St}(\mathfrak A) \ar[r]&
%E(\mathfrak A) \ar[r]& (1)}
%\tag{5}
%\end{equation*}

\begin{equation}
\begin{CD}1@>>>$K$_2(\mathcal A)@>>>  {\rm St}(\mathcal A)@>>>E(\mathcal A)@>>>1\end{CD}
\tag{10}
\end{equation}
o\`u $E(\mathcal A)$ d\'esigne le sous-groupe du groupe lin\'eaire infini
engendr\'e par les matrices \'el\'ementaires, ${\rm St}(\mathcal A)$ son extension centrale universelle ou groupe de Steinberg, et $K_2(\mathcal A)$ le deuxi\`eme groupe de $K$-th\'eorie de $\mathcal A$ (voir Milnor \cite{Milnor}).

Nous voici, du moins en apparence, \'eloign\'es du cocycle de
Virasoro. Revenons \`a une alg\`ebre de Lie $\g$ munie d'une forme
invariante $\kappa : \g \times \g \to \mathbb{K}$, et soit $\delta : \g \to \g$
une d\'erivation ext\'erieure. Il est facile de v\'erifier que si
l'application $(X,Y) \mapsto \kappa (X, \delta Y)$ est
antisym\'etrique, alors elle d\'efinit un $2$-cocycle que l'on notera
$c_\delta$. En particulier si $\g$ est une alg\`ebre de Lie obtenue
par antisym\'etrisation d'une alg\`ebre associative, l'existence d'une
trace permet d'obtenir une forme bilin\'eaire invariante suivant la formule:
$\kappa (X,Y)={\rm Tr}(XY)$. C'est le cas des alg\`ebres de Kac-Moody
en consid\'erant chaque alg\`ebre de Lie semi-simple de dimension
finie comme sous-alg\`ebre d'un $\mathfrak{gl}(n)$; pour les
alg\`ebres de courants $\g_{\mathcal A}$, on peut obtenir une trace en
utilisant un ``r\'esidu'' ou une ``int\'egrale'' : $\Omega^1 \mathcal A
\to \mathbb{K}$ ; enfin, rappelons l'existence de la{\it trace d'Adler }
sur l'alg\`ebre associative $\psi \mathcal D (S^1)$ des symboles
d'op\'erateurs pseudo-diff\'erentiels sur $S^1$, 
$$
{\rm Tr}\Big(\sum_{i=-\infty}^n a_i (x)\xi^i\Big)=\int_{S^1} a_{-1}(x)dx ,
$$
qui permet de construire deux extensions centrales
ind\'ependantes avec les d\'erivations ext\'erieures ${\rm ad \,
Log}\xi$ et ${\rm ad}\, x$ qui donnent les cocycles  $c_1 (D_1,
D_2)={\rm Tr}(D_1 [{\rm Log}\xi,D_2])$ et $c_2 (D_1,D_2)={\rm
Tr}(D_1[x,D_2])$ (r\'esultat de O. Kravchenko et B. Khesin). Rappelons
que le plongement naturel ${\rm Vect}(S^1) \to \psi \mathcal D (S^1)$
permet d'induire de $c_1$ le cocycle de Virasoro.

En fait, tous ces multiples cocycles \`a valeurs scalaires, dont ceux de Virasoro et Kac-Moody,  peuvent
s'obtenir par plongement de l'alg\`ebre de Lie concern\'ee dans
{\it l'alg\`ebre des matrices de Jacobi} et restriction de son cocycle
universel. Rappelons bri\`evement la formule de ce dernier: les matrices de
Jacobi sont des matrices doublement infinies $a_{ij}$ pour $(i,j) \in
\Z^2$ qui v\'erifient $a_{ij}=0$ si $\abs{i-j}$ est assez grand; il
est facile de v\'erifier que ces matrices forment une alg\`ebre
associative et donc une alg\`ebre de Lie par antisym\'etrisation. Si
l'on note alors $J$ la matrice d\'efinie par $a_{ij}=\delta_{ij}{\rm
sgn}(i)$ pour $j \ne 0$ et $a_{00}=1$, le cocycle va s'\'ecrire
\begin{equation}
\label{11}
\displaystyle c(A,B)=\frac{1}{2}{\rm Tr}(A[J,B])
\tag{11}
\end{equation}

Une telle formule est a priori d\'epourvue de sens car la trace n'est
pas d\'efinie pour une matrice de Jacobi, mais ici un calcul direct
montre que $A[J,B]$ est en fait une matrice finie (tout le monde aura
compris qu'il s'agit d'une matrice dont tous les \'el\'ements sont
nuls sauf un nombre fini d'entre eux). On obtient une formule
explicite  en d\'ecomposant chaque matrice de Jacobi en blocs
$\begin{bmatrix}A_{11} & A_{12} \cr A_{21} & A_{22} \end{bmatrix}$; on obtient
$c(A,B)={\rm Tr}(A_{21} B_{12} -A_{12}B_{21})$, et la
v\'erification de la propri\'et\'e de cocycle devient un exercice
\'el\'ementaire.

Nous avons ainsi construit une extension centrale non triviale 
%\begin{equation*}
%\label{7}
%\xymatrix{0 \ar[r]& \mathbb{K} \ar[r]& \wh{\mathfrak{gl}_J}(\mathfrak A) \ar[r]& 
%\mathfrak{gl}_J (\mathfrak A) \ar[r]& 0 }
%\tag{7}
%\end{equation*} 
\begin{equation}
\displaystyle
\begin{CD}
0@>>> \mathbb{K} @>>>\hat{\mathfrak{gl}_J}(\mathcal A)@>>> \mathfrak{gl}_J (\mathcal A)@>>> 0\end{CD}
\tag{12}
\end{equation}

Mais nous n'avons pas pr\'ecis\'e quel doit \^etre l'anneau de base
$\mathcal A$ ; il suffit que $\mathcal A$ soit une alg\`ebre commutative sur
le corps $\mathbb{K}$ avec une int\'egrale ou un r\'esidu permettant de
d\'efinir une trace \`a valeurs scalaires 
${\rm Tr} : \mathfrak{gl}_{\rm fini} (\mathcal A)\to \mathbb{K}$. 

Une extension analogue
peut se construire dans une version analytique, avec des espaces de Hilbert,
ce qui permet l'int\'egration de l'extension de l'alg\`ebre au groupe lin\'eaire.
 On part d'un espace de Hilbert d\'ecompos\'e suivant
$\mathcal H=\mathcal H_+ \oplus \mathcal H_-$, et en posant $J_{|\mathcal
H_\pm}=\pm {\rm Id}_{\mathcal H_\pm}$; on d\'efinit ensuite le groupe lin\'eaire
restreint ${\rm GL}_{\rm res} (\mathcal H)$ comme le sous-groupe des
applications lin\'eaires born\'ees $A : \mathcal H \to \mathcal H$ telles
que $[J,A]$ soit de Hilbert-Schmidt; on remarque alors que pour $A$
et $B$ dans ${\rm GL}_{\rm res}(\mathcal H)$, $A[J,B]$ est tra\c cable;
autrement dit ${\rm Tr}(A [J,B])$ est bien d\'efinie. On en d\'eduit
les extensions au niveau du groupe et de l'alg\`ebre de Lie 

%\begin{equation*}
%\xymatrix{(0) \ar[r]& \C^* \ar[r]& \wh{{\rm GL}_{\rm res}}(\mathfrak H)
%\ar[r]& {\rm GL}_{\rm res} (\mathfrak H) \ar[r]& (0) }
%\tag{8}
%\end{equation*}
\begin{equation}
\displaystyle
\begin{CD}
1 @>>>\mathbb{C}^*@>>> \hat{\rm GL}_{\rm res}(\mathcal H)@>>>{\rm GL}_{\rm res} (\mathcal H) @>>> 1\end{CD}
\tag{13}
\end{equation}
%\begin{equation*}
%\label{9}
%\xymatrix{(0) \ar[r]& \C \ar[r]& \wh{\mathfrak{gl}_{\rm res}}(\mathfrak H)
%\ar[r]& \mathfrak{gl}_{\rm res} (\mathfrak H) \ar[r]& (0) }
%\tag{9}
%\end{equation*}
\begin{equation}
\displaystyle
\begin{CD}
0 @>>>\mathbb{C}@>>>  \hat{\mathfrak{gl}_{\rm res}}(\mathcal H)@>>>\mathfrak{gl}_{\rm res} (\mathcal H) @>>> 0\end{CD}
\tag{14}
\end{equation}
L'extension (13) est en quelque sorte la compl\'etion de
(11). Le cocycle de Bott-Thurston et celui de Virasoro pour ${\rm
Diff}(S^1)$ peuvent \^etre obtenus via un plongement de ${\rm
Diff}(S^1)$ dans ${\rm GL}_{\rm res}(\mathcal H)$, respectivement de ${\rm
Vect}(S^1)$ dans $\mathfrak{gl}_{\rm res}(\mathcal H)$, d\'eduit de
l'action de ${\rm Diff}(S^1)$ sur $\mathcal H=L^2 (S^1, \C)$ par
changement de variables. Ce cocycle a une parfaite
analogie formelle avec celui de (4), il est de la forme
$c(A,B)={\rm Tr}(A \delta B)$ o\`u $\delta$ est une d\'erivation;
cependant cette analogie n'est que formelle, car la trace n'est pas
d\'efinie globalement sur $\mathfrak{gl}_J$ ou $\mathfrak{gl}_{\rm
res}$ et l'application $A \mapsto [J,A]$ n'est pas une d\'erivation
ext\'erieure. Ce cocycle a \'et\'e g\'en\'eralis\'e par Alain Connes
qui en a fait l'une des pierres angulaires de sa g\'eom\'etrie non
commutative \cite{Connes1,Connes2} dans le cadre des modules de Fredholm,
qui sont des espaces de Hilbert gradu\'es $\mathcal H=\mathcal H_+ \oplus \mathcal H_-$.
Une formule $\displaystyle {\rm Tr}_S (A)=\frac{1}{2}{\rm Tr}(\gamma
F[F,A])$ permet de d\'efinir une trace, o\`u $\gamma$ est un
op\'erateur de graduation sur $\mathcal H$. Le cocycle pr\'ec\'edent peut
alors s'\'ecrire sous la forme $c(A,B)={\rm Tr}_S (A[F,B])$; il
appara\^\i t comme un exemple de classe de cohomologie cyclique
\cite{Loday1}, et en quelque sorte le prototype de celles-ci:
dans\cite{Connes1}, p 196, la formule d'un cocycle cyclique associ\'e \`a un module de Fredholm
s'\'ecrit $c(A_0,A_1,....,A_n)={\rm Tr}_S (A_0[F,A_1]\ldots[F,A_n])$. L'op\'erateur $A
\mapsto dA=[F,A]$ porte le nom de diff\'erentielle quantique ; l'analogie avec le cocycle (4) est
ainsi rendue plus pr\'ecise. 

\smallskip

Le surnom de ``cocycle japonais'', donn\'e parfois \`a ce cocycle par l'�\'ecole
russe"(celle de St-Petersbourg, cf. plus haut), provient de l'intervention d\'ecisive de ce cocycle dans les travaux de
l'\'ecole dite de Kyoto; les math\'ematiciens de cette ville, travaillant au RIMS autour du Pr Mikio Sato c\'el\`ebre pour ses formules, sont bien connus,
M. Jimbo et ses collaborateurs notamment, pour leurs r\'esultats sur les syst\`emes
int\'egrables en dimension infinie, et notamment les hi\'erarchies KdV et KP. La repr\'esentation d'une alg\`ebre de Lie de dimension infinie judicieusement choisie permet d'engendrer des familles de solutions de ces \'equations \`a solitons, et le cocycle y joue un r\^ole de 'deus ex machina' \cite{Dickey,KaRa}. Cependant
on ne peut  attribuer \`a cette "\'ecole de Kyoto" la priorit\'e dans la d\'ecouverte de ce cocycle; sa
premi\`ere description dans un cadre math\'ematique pr\'ecis se trouve dans le
c\'el\`ebre trait\'e de Pressley et Segal \cite{PressleySegal}. 
Tel ou tel
math\'ematicien pourra en revendiquer la paternit\'e, mais sa
v\'eritable origine se trouve en physique dans les travaux de
F. Berezin (\cite{Berezin}, 1965) sur la seconde quantification( voir\cite{Roger} pour un survol des multiples contributions de F. Berezin \`a la physique math\'ematique). 
Y. Neretin \cite{Neretin}
affirme toutefois que la premi\`ere formule explicite pour ce cocycle
est due \`a K.O. Friedrichs \cite{Friedrichs}, dans son \'etude du " groupe canonique" dont il va \^etre question maintenant... 

\smallskip

Ceci nous fournit une transition
naturelle vers les origines physiques du cocycle de Virasoro. Le passage de la physique classique \`a la physique quantique se fait d'abord de la m\' ecanique classique \`a la m\'ecanique quantique,  c'est \`a dire dans le cas d'un nombre fini de degr\'es de libert\'e; l'espace des \' etats quantiques se pr\'esente alors sous la forme d'un espace de Hilbert $\mathcal H$, celui des fonctions d'onde. La th\'eorie quantique des champs concerne un nombre infini de degr\'es de libert\'e, et on sait qu'il n'en existe pas encore de th\'eorie math\'ematique naturelle et pleinement satisfaisante ( pour un expos\'e lisible et pragmatique, voir\cite{Zee}). Une premi\`ere et rudimentaire version en est la th\'eorie du champ libre, que l'on peut voir comme la repr\'esentation quantique d'une assembl\'ee de particules, ou "seconde quantification". Rappelons maintenant le principe de cette seconde quantification: on associe \`a l' espace de Hilbert des \'etats $\mathcal H$, un {\it espace de Fock} $\mathcal F(\mathcal H)$ (fermionique ou
bosonique suivant le comportement statistique des particules) et on essaie de prolonger \`a $\mathcal F(\mathcal H)$ l'action d'un groupe de
sym\'etries de $\mathcal H$.

Pr\'ecisons tout d'abord ce qu'est l'espace de Fock, dans le cas bosonique: soit $\mathcal H$ est un espace hermitien de dimension $n$, munissons-le de
la mesure gaussienne $\displaystyle d \mu (z)=\pi^{-n}\exp(-\norm{z}^2)\, dz$ ; l'espace de Fock bosonique $\mathcal F(\mathcal H)$ est alors l'espace des
fonctions holomorphes sur $\mathcal H$ de carr\'e sommable pour cette
mesure. C'est un espace de Hilbert et toute base de $\mathcal H$ permet
de d\'efinir une base de $\mathcal F(\mathcal H)$: les fonctions
coordonn\'ees $z_1, \ldots, z_n$ associ\'ees \`a une base $\{e_1,
\ldots, e_n\}$ de $\mathcal H$ d\'efinissent des fonctions  pour tous les
multi-indices $I \in \mathbb {N}^n$ par la formule 
$$
e_I (z) = z^{i_{1}}_1 z^{i_{2}}_2 \dots z^{i_{n}}_n 
\qquad
\text{ si }I =(i_{1},i_{2},\dots,i_{n}).
$$
Les $e_I (z)$ constituent une base
orthonorm\'ee de $\mathcal F(\mathcal H)$ comme on le v\'erifie facilement
par un calcul d'int\'egrales gaussiennes ; il est alors clair que
l'anneau des polyn\^omes $\mathbb{ C}[z_1, \ldots , z_n]$ se plonge dans
$\mathcal F(\mathcal H)$ comme sous-espace partout dense.

On peut ensuite \'etendre cette construction de l'espace de Fock $\mathcal F(\mathcal H)$ \`a tout
espace de Hilbert complexe s\'eparable, par un proc\'ed\'e de passage \`a la limite, ou alors associer \`a une base hilbertienne $\{e_i\}_{i \in \mathbb{N}}$ de $\mathcal H$,  une base $e_I$ de la limite inductive des
$\mathcal F({\mathcal H})$(o\`u $\mathcal H$ est de dimension finie comme
ci-dessus), et d\'efinir $\mathcal F(\mathcal H)$ comme l'espace de Hilbert (non s\'eparable) des
$\displaystyle \sum_{I \in \mathbb{N}^{\infty}}a_I e_I$ avec $\sum \abs {a_I}^2<+\infty$.  Une autre option, peut-\^etre plus intuitive consiste \`a proc\'eder alg\' ebriquement, en posant
$\mathcal F(\mathcal H)=\underset{k=0}{\overset{+
\infty}{\oplus}}{\cal S}^k(\mathcal H)$ o\`u ${\cal S}^k(\mathcal H)$ d\'esigne l'espace des
$k$-tenseurs sym\'etriques sur $\mathcal H$ , puis \`a compl\'eter cet espace avec les pr\'ecautions d'usage.

Si maintenant on essaie de prolonger l'action d'un groupe de
sym\'etries quantiques de $\mathcal H$ \`a $\mathcal F(\mathcal H)$, on peut trouver une solution \`a condition de se limiter au groupe lin\'eaire restreint (symplectique ou orthogonal
suivant le cas) et la repr\'esentation sur $\mathcal F (\mathcal H)$ est
projective et non plus lin\'eaire ; et c'est ainsi que le cocycle
(6) est apparu naturellement, comme terme central d'une extension de ce ``{\it groupe
canonique}'', comme on dit parfois dans la litt\'erature
physique. Voici pr\'ecis\'ement venu le moment de parler des
``{\it anomalies}'', qui ne sont autres que les classes de cohomologie qui
appara\^{\i}ssent dans certains calculs de la th\'eorie quantique
des champs...

\smallskip 

Une anomalie pourra se pr\'esenter comme un d\'efaut \`a \^etre un homomorphisme, ou
un d\'efaut \`a \^etre une action de groupe, ou encore une
courbure... Les dites anomalies se repr\'esentent par des classes de
cohomologie appropri\'ees, mais l'\'elaboration du cadre n\'ecessaire
pour ces cohomologies (cohomologie de groupes, d'alg\`ebres de Lie, ou
th\'eories topologiques...)est souvent fort d\'elicate, et a parfois
suscit\'e des probl\`emes math\'ematiques originaux et
difficiles... Ces anomalies pourront \^etre ``globales'',
c'est-\`a-dire repr\'esenter des invariants topologiques, du type
classes caract\'eristiques ou invariants de Chern-Simons, comme l'action de Wess-Zumino-Witten, ou alors ``locales'', purement
alg\'ebriques, et se ramenant \`a des classes de cohomologie de
certains groupes ou alg\`ebres de Lie de dimension infinie. Du point de vue
lagrangien, la distinction local-global appara\^{\i}t suivant que les termes
suppl\'ementaires modifient ou non la dynamique. Ces d\'efauts appara\^{\i}ssent lors
du passage classique
$\to$ quantique et sont intrins\`equement li\'es \`a la quantification. Une situation typique est la suivante : au niveau classique on a
un homomorphisme $\phi$ d'alg\`ebres de Poisson, donc :
$$
\{\phi (f), \phi (g)\}=\phi (\{f,g\})
$$
et sa quantification $\hat{\phi}$ pr\'esente un d\'efaut \`a \^etre un
homomorphisme suivant la formule:
$$
[\hat{\phi}(f),\hat{\phi}(g)]=\hat{\phi}([f,g])+{\rm Id}\, \, c(f,g),
$$
faisant ainsi appara\^\i tre un cocycle $c$ \`a valeurs
scalaires.

 Retracer l'histoire des anomalies pourrait faire l'objet d'une monographie enti\`ere; nous en indiquerons seulement les principales
\'etapes, en d\'eveloppant les points en relation directe avec le
cocycle de Virasoro. Ce fut d\`es les ann\'ees 50 que Julian Schwinger
mit en \'evidence un ph\'enom\`ene d'anomalies, dans le cadre de
l'\'electrodynamique quantique, pour l'\'etude de la
d\'esint\'egration du $\pi^0$ (m\'eson-pi ou pion) en photons ;
%%%%%%%%%%%%%%%%%%%%%%%%
\begin{figure}[htbp]
\begin{center}
 \begin{fmfgraph*}(100,100)
     \fmfpen{thick}
      \fmfleft{i}
      \fmfright{o1,o2}
      \fmf{plain,label=$\pi^0$,label.side=left}{i,v}
      \fmfpen{thin}
      \fmf{photon,label=$\gamma$,label.side=right}{v,o1}
     \fmf{photon,label=$\gamma$,label.side=left}{v,o2}
     	%\fmffreeze
      %\fmf{fermion,width=.6thick}{v,o3}
      %\fmf{photon}{v3,v2}
   %   \fmffreeze
      %\fmf{photon,tension=0}{v3,o2}
%      \fmflabel{$\pi$}{i1}
   \end{fmfgraph*}
   \caption{D\'esint\'egration du pion neutre}
\label{fig:mongraphe}
\end{center}
\end{figure}
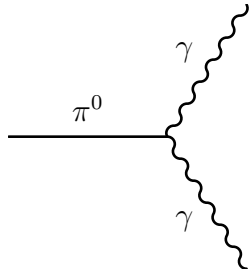
le calcul des propagateurs suivant le diagramme de Feynman ci-dessous fait appara\^\i tre une anomalie, c'est \`a dire une obstruction \`a \'etendre  les sym\'etries classiques au cas quantique. Le calcul est d\'etaill\'e dans le trait\'e de th\'eorie quantique des champs de Peskin et Schroeder, on y trouve le r\'esultat final 
dans la formule 19.45  \cite{Peskin}[p 661]; cette formule peut s'\'ecrire en langage g\'eom\'etrique sous la forme:

$ Div(J)=\frac {e^2}{16\pi^2}(<\Lambda,\Omega\wedge\Omega>)$

o\`u $\Lambda$ est le tenseur de volume contravariant,
$\Omega$ d\'esigne l'intensit\'e du champ \'electromagn\'etique, et bien entendu, $J$ est le vecteur courant. Le calcul quantique met donc en d\'efaut la conservation du courant, et c'est l\`a que se trouve l'anomalie!
%%%%%%%%%
On peut interpr\'eter math\'ematiquement cette anomalie comme une singularit\'e dans un calcul de
distributions, soit un produit de deux distributions de Dirac, ou encore comme un
d\'efaut d'homomorphisme mettant en \'evidence un cocycle central, le
groupe sous-jacent \'etant le groupe de jauge de l'\'electrodynamique,
soit $C^\infty (X, {\rm U}(1))$ pour un espace-temps $X$\cite{Schwinger2,Schwinger3,Schwinger4}. L'\'etude de
ces ph\'enom\`enes fut ensuite poursuivie par R. Jackiw et S. Adler,
cr\'eateurs du terme d'anomalies \cite{AdBa} en 1969, toujours pour
l'\'electrodynamique quantique. Au m\^eme moment, on assistait \`a un d\'eveloppement spectaculaire
des th\'eories de jauge non ab\'eliennes (Yang-Mills) avec entre autres les courants faibles et la chromodynamique quantique; le ph\'enom\`ene d'anomalies dans
le cadre non ab\'elien fut mis en \'evidence par S. Adler et
D. Bardeen; c'est aussi vers cette \'epoque qu'au cours de l'un des
multiples colloques consacr\'es \`a la th\'eorie des groupes en
physique, Andr\'e Lichn\'erowicz pouvait proph\'etiser que la physique
th\'eorique allait devenir de plus en plus cohomologique....

 La th\'eorie g\'en\'erale des anomalies, et sa compr\'ehension dans le cadre alg\'ebrique et g\'eom\'etrique ad\'equat, fut \'elabor\'ee presque simultan\'ement
et ind\'ependamment par L.D. Fadeev, R. Stora et B. Zumino (ordre
alphab\'etique). Du point de vue math\'ematique, il s'agit de classes
de cohomologie du groupe de jauge $C^\infty (V,G)=\mathcal G$ ou de l'alg\`ebre des courants
$C^\infty (V ; \g)=\underline{\mathcal G}$,  o\`u $G$ est un groupe compact non
ab\'elien. Les articles \cite{St,Fa} constituent une approche
int\'eressante et accessible pour un math\'ematicien; par
exemple on voit dans \cite{Fa} appara\^\i tre le cocycle comme une obstruction
dans l'\'ecriture de la quantification canonique d'une th\'eorie de
jauge. M\^eme si ces cocycles ne sont pas directement li\'es \`a celui de Virasoro sauf par l'interm\'ediaire de leur "container" commun longuement d\'ecrit plus haut, leur contexte et leur d\'etermination sont tout \`a fait analogues, et nous allons donc en donner certains d\'etails.

Les formules de \cite{Fa} peuvent para\^\i tre un peu myst\'erieuses
au premier abord, mais c'est un exercice instructif que de les transcrire
sous une forme alg\'ebrico-g\'eom\'etrique plus globale. Prenons par
exemple la formule (22)
$$[G^a (x), G^b (y)]=if^{abc} G^c (x) \delta^{(3)} (x-y)+\frac{1}{12 i
\pi^2} d^{abc}\varepsilon_{ijk}\partial_i A^c_j (x) \partial_k
(\delta^{(3)}(x-y))$$
o\`u $f^{abc}$ d\'esigne les constantes de structure de l'alg\`ebre de
Lie simple $\mathfrak g$, $d^{abc}$ d\'esigne les coefficients du g\'en\'erateur de $H^3
(\g, \C)$, $\underline{A}=A_j^c (x)e_c dx_j$
d\'esigne le champ de jauge, soit en langage math\'ematique  la forme
de connexion. L'espace-temps est ici de dimension $3$ et
$\varepsilon_{ijk}$ d\'esigne sa forme volume. L'anomalie est alors
le second terme dans le second membre. Dans une \'ecriture globale,
cette formule nous donne: pour $f$ et $g$ dans $\mathcal G$ on a 
\begin{equation*}
\label{}
\displaystyle [f,g](x)=[f(x), g(x)]+\frac{1}{12 i \pi^2}\int_V T (fdg
\wedge d \underline{A})
\end{equation*}
o\`u la forme $T(fdg \wedge \underline{A})$ s'interpr\`ete comme 
$$T(fdg \wedge \underline{A})(x)=d^{abc} f^a (x)dg^b (x) \wedge dA^c
(x).$$
En notant $\mathcal C$ l'espace des champs de jauge et $\mathcal F$
l'espace des fonctions sur $\mathcal C$, on peut consid\'erer cette
anomalie comme un $2$-cocycle \`a valeurs dans $\mathcal F$
\begin{equation*}
\begin{aligned}
\underline{\mathcal G} \times \underline{\mathcal G} & \longrightarrow \mathcal
F\\
(f,g) & \longmapsto c (f,g)
\end{aligned}
\end{equation*}
o\`u $c(f,g)$ est d\'efini par $\displaystyle
c(f,g)(\underline{A})=\int_V T (fdg \wedge \underline{A})$.

\smallskip 

Ces r\'esultats sont curieusement analogues et presque contemporains
des calculs cohomologiques sur $\g_{\mathcal A}$ que nous avons mentionn\'es plus haut
; ils ont \'et\'e en tout cas  \'elabor\'es de fa\c con strictement ind\'ependante et l'analogie ne semble pas
avoir
\'et\'e remarqu\'ee \`a l'\'epoque. 

\smallskip

Un peu plus tard, Mickelsson a donn\'e de ce cocycle
une pr\'esentation plus analytique, en terme de seconde
quantification, permettant de trouver l'anomalie \`a partir du cocycle
universel du groupe lin\'eaire restreint ; pour chaque connexion $A
\in \mathcal C$, on a un espace de Fock fermionique $\mathcal F_A$
associ\'e, construit \`a partir des solutions de l'\'equation de Dirac
pour $A$ ; les transformations de jauge agissent alors simultan\'ement sur
la connexion et sur l'espace de Fock et c'est ce qui fait toute la
difficult\'e de ces th\'eories de jauge. On en d\'eduit une
repr\'esentation $\underline{\mathcal G} \to \mathfrak{gl}_{\rm
res}(\mathcal H)$ et l'image r\'eciproque du cocycle (9) sur 
$\mathfrak{gl}_{\rm res}$ nous donne le cocycle de $\underline{\mathcal
G}$. Ce cocycle est connu dans la litt\'erature sous le nom de cocycle
de Mickelsson-Rajeev \cite{MiRa} lorsque la dimension de l'espace-temps est
\'egale \`a $3$. C'est un habillage analytique du cocycle de Faddeev
que nous avons d\'ecrit plus haut. Il peut  en outre s'exprimer en terme d'op\'erateurs
pseudo-diff\'erentiels et du r\'esidu de Wodzicki. Un expos\'e tr\`es
clair de tous ces r\'esultats se trouve dans les travaux de
C. Ekstrand \cite{Ek}. 

\bigskip

Dans le m\^eme ordre d'id\'ees, nous allons maintenant mentionner bri\`evement un
travail de G. Segal intitul\'e ``Faddeev's anomaly in Gauss's law'',
tr\`es int\'eressant mais h\'elas non publi\'e; G. Segal y donne
l'analyse topologique de ces anomalies en th\'eorie de jauge. Les
espaces de Fock fermioniques forment un fibr\'e sur l'espace des
connexions $\mathcal C$, soit $\displaystyle \mathcal F
\overset{\pi}{\longrightarrow} \mathcal C$ avec $\pi^{-1}(A)=\mathcal
F_A$. Cet espace $\mathcal F_A$ est isomorphe naturellement \`a $\mathcal
F$ c'est l'action de jauge qui est donn\'ee via le potentiel $A$.Le groupe de jauge op\`ere donc naturellement et on a une projection
des espaces quotients $\mathcal F/\mathcal G \overset{\pi}{\longrightarrow}
\mathcal C /\mathcal G$. L'espace des \'etats physiques est alors l'espace
des sections de ce ``fibr\'e'' ; ce n'est pas un v\'eritable fibr\'e
en vari\'et\'es de dimension infinie \`a cause des singularit\'es de
l'action de $\mathcal G$ sur $\mathcal C$ (voir par exemple les travaux de
Singer \cite{Si}). G. Segal d\'eduit ensuite du projectifi\'e $\mathcal P$
du fibr\'e de Fock $\mathcal F$, un fibr\'e $\mathcal P/\mathcal
G \overset{\pi}{\longrightarrow} \mathcal C /\mathcal G$ de fibre $P(\mathcal
H)$ o\`u $\mathcal H$ est un espace de Hilbert complexe; cet espace
$P(\mathcal H)$ est isomorphe \`a un projectif complexe infini $\C
P(\infty)$, et c'est donc topologiquement un espace d'Eilenberg-Mac Lane
$K(\Z, 2)$; la projection $\mathcal P/\mathcal G
\overset{\pi}{\longrightarrow} \mathcal C /\mathcal G$ est donc
caract\'eris\'ee par une classe de cohomologie enti\`ere $\widetilde{c} \in H^3_{\rm
top} (\mathcal C/\mathcal G, \Z)$. L'auteur montre ensuite par un argument
tr\`es subtil, du type th\'eor\`eme de Van-Est, la relation entre
cette classe topologique et le cocycle de Faddeev: on a une
application $ H^3_{\rm top} (\mathcal C/\mathcal G,
\Z)\overset{I}{\longrightarrow} H^2 (\underline{\mathcal G}, \mathcal F
(\mathcal C, \R))$ d\'efinie naturellement; soit une forme ferm\'ee  $\omega \in \Omega^3
(\mathcal C/\mathcal G)$ repr\'esentant $\widetilde{c}$,
elle se rel\`eve suivant une forme ferm\'ee $\overline{\omega} \in \Omega^3
(\mathcal C)$. L'espace $\mathcal C$ \'etant convexe, donc acyclique, il
existe $\alpha \in \Omega^2 (\mathcal C)$ telle que $d \alpha
=\overline{\omega}$; la classe $I(\widetilde{c})$ va alors se
d\'efinir  par

$$I(\widetilde{c})(X,Y)(A)=\overline{\omega}_A (\theta_X (A),\theta_Y (A))$$

o\`u $\theta_X (A)$ d\'esigne l'action infinit\'esimale du courant $X$
sur la connexion $A$ ; en d'autres termes, $\theta_X (A)=dX +
[X,A]$. G. Segal montre enfin que $I(\widetilde{c})$ s'identifie au
cocycle de Faddeev. Nous nous sommes quelque peu attard\'e sur ce
sujet m\^eme s'il nous \'eloigne en apparence du cocycle de Virasoro, car il est tr\`es instructif pour comprendre alg\'ebriquement les anomalies. 
\bigskip

C'est maintenant qu'il faut parler des {\it mod\`eles
duaux}, assez oubli\'es de nos jours, mais qui eurent leur heure de
gloire dans les ann\'ees 70 et se d\'evelopp\`erent simultan\'ement
et en interaction (si on ose dire) avec les diverses anomalies.

\begin{figure}[htbp]
\begin{center}
\subfloat[ $1$]{
\begin{fmfgraph*}(100,100)
      \fmfleft{i1,i2,i3}
      \fmfright{o1,o2,o3}
      \fmfpen{thin}
      \fmf{dashes,width=thin}{i2,v1}
      \fmf{dashes,tension=.1,width=thin}{v1,v2}
            \fmf{dashes,width=thin}{v2,o2}
        	\fmf{plain,thick,left,tension=0.2,tag=1}{v1,v2} 
	\fmf{plain,left,tension=0.2,tag=2}{v2,v1} 
	
	\fmfforce{c+(.15w*(cosd(45),sind(45)))}{w1}
	\fmfforce{c+(.15w*(cosd(135),sind(135)))}{w2}
	\fmfforce{c+(.15w*(cosd(-45),sind(-45)))}{w3}
	\fmfforce{c+(.15w*(cosd(-135),sind(-135)))}{w4}
	%\fmfdotn{w}{4} 	
	\fmf{fermion,width=.6thick}{i1,w4}
	\fmf{fermion,width=.6thick}{i3,w2}
	\fmf{fermion,width=.6thick}{o3,w1}
	\fmf{fermion,width=.6thick}{o1,w3}
	\fmflabel{$\mathbf{P_1}$}{i1}
	\fmflabel{$\mathbf{P_2}$}{i3}
	\fmflabel{$\mathbf{P_3}$}{o3}
	\fmflabel{$\mathbf{P_4}$}{o1}
    \end{fmfgraph*}\label{fig:ex1}}\hspace{2cm}%
    \subfloat[ $2$]{
\begin{fmfgraph*}(100,100)
      \fmfleft{i1,i2,i3}
      \fmfright{o1,o2,o3}
      \fmftop{t}
      \fmfbottom{b}
      \fmfpen{thin}
      
      \fmf{dashes,width=thin}{t,v1}
          \fmf{dashes,tension=.1,width=thin}{v1,v2}
         \fmf{dashes,width=thin}{v2,b}
        	
	\fmf{plain,thick,left,tension=0.2,tag=1}{v1,v2} 
	\fmf{plain,left,tension=0.2,tag=2}{v2,v1} 
	
	\fmfforce{c+(.15w*(cosd(45),sind(45)))}{w1}
	\fmfforce{c+(.15w*(cosd(135),sind(135)))}{w2}
	\fmfforce{c+(.15w*(cosd(-45),sind(-45)))}{w3}
	\fmfforce{c+(.15w*(cosd(-135),sind(-135)))}{w4}
	%\fmfdotn{w}{4} 	
	\fmf{fermion,width=.6thick}{i1,w4}
	\fmf{fermion,width=.6thick}{i3,w2}
	\fmf{fermion,width=.6thick}{o3,w1}
	\fmf{fermion,width=.6thick}{o1,w3}
	\fmflabel{$\mathbf{P_1}$}{i1}
	\fmflabel{$\mathbf{P_2}$}{i3}
	\fmflabel{$\mathbf{P_3}$}{o3}
	\fmflabel{$\mathbf{P_4}$}{o1}
    \end{fmfgraph*}\label{fig:ex2}}  
    \caption{Mod\`ele duaux}
\label{fig:mongraphe}
\end{center}
\end{figure}
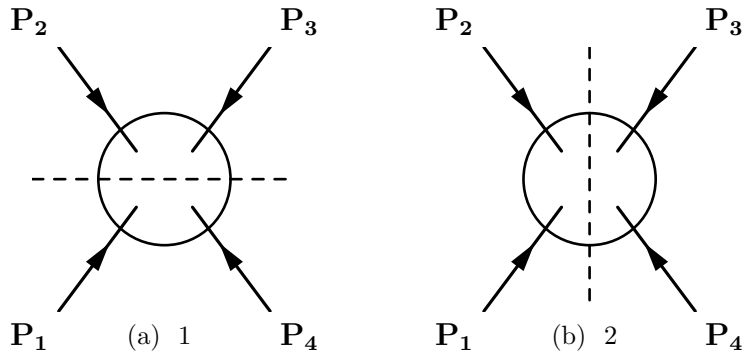

Indiquons bri\`evement l'id\'ee de d\'epart de cette th\'eorie des
mod\`eles duaux, qui a \'et\'e le premier lieu de l'apparition de l'alg\`ebre de
Virasoro en physique. Dans l'\'etude des interactions de deux
particules lourdes (les {\it hadrons}) repr\'esent\'ee par un diagramme du
type suivant (\ref{fig:ex1}) soit $P_1+P_2 \to P_3+P_4$; l'id\'ee de base des
mod\`eles duaux consiste \`a consid\'erer simultan\'ement ce m\^eme diagramme
comme celui de l'interaction $P_2 +P_3 \to P_1 + P_4$, soit (\ref{fig:ex2}).

%\mydessf{20}{60}{dua3.eps}{0}{dua}{}
%\begin{minipage}[.3\textwidth]
%\end{minipage}

%%\begin{minipage}[.3\textwidth]
%\begin{figure}[htbp]
%\begin{center}
%\begin{fmfgraph*}(100,100)
%      \fmfleft{i1,i2,i3}
%      \fmfright{o1,o2,o3}
%      \fmfpen{thin}
%      \fmf{dashes,width=thin}{i2,v1}
%      \fmf{dashes,tension=.1,width=thin}{v1,v2}
%            \fmf{dashes,width=thin}{v2,o2}
%        	\fmf{plain,thick,left,tension=0.2,tag=1}{v1,v2} 
%	\fmf{plain,left,tension=0.2,tag=2}{v2,v1} 
%	
%	\fmfforce{c+(.15w*(cosd(45),sind(45)))}{w1}
%	\fmfforce{c+(.15w*(cosd(135),sind(135)))}{w2}
%	\fmfforce{c+(.15w*(cosd(-45),sind(-45)))}{w3}
%	\fmfforce{c+(.15w*(cosd(-135),sind(-135)))}{w4}
%	%\fmfdotn{w}{4} 	
%	\fmf{fermion,width=.6thick}{i1,w4}
%	\fmf{fermion,width=.6thick}{i3,w2}
%	\fmf{fermion,width=.6thick}{o3,w1}
%	\fmf{fermion,width=.6thick}{o1,w3}
%	\fmflabel{$\mathbf{P_1}$}{i1}
%	\fmflabel{$\mathbf{P_2}$}{i3}
%	\fmflabel{$\mathbf{P_3}$}{o3}
%	\fmflabel{$\mathbf{P_4}$}{o1}
%    \end{fmfgraph*}
% \caption{mod\`ele dual (suite)}
%\label{fig:mongraphe}
%\end{center}
%\end{figure}
%%\end{minipage}

La loi de conservation de l'impulsion donne $P_1+P_2+P_3+P_4=0$ et on
en d\'eduit les variables d'\'energie (``de Mandelstam'')
\begin{xalignat*}{2}
S  =(P_1+P_2)^2=(P_3+P_4)^2 & \qquad t  =(P_2 +P_3)^2=(P_1+P_4)^2.
\end{xalignat*}
La fonction d'amplitude $A(S,t)$ est alors une fonction analytique de
$S$ et de $t$ ; on voit appara\^\i tre des ph\'enom\`enes de
r\'esonance d'o\`u le nom de \textit{dual resonance models} pour ces
th\'eories. Du point de vue ph\'enom\'enologique, ces r\'esonances
appara\^\i ssaient comme de nouvelles particules hadroniques, et leur
nombre semblait tr\`es grand; le d\'eveloppement de la physique exp\'erimentale des particules au cours des ann\'ees soixante a permis la d\'ecouverte d'un nombre imposant de ces particules lourdes, que la th\'eorie alors disponible ne permettait pas encore de classifier et d'interpr\`eter de mani\`ere satisfaisante; pour un historique du sujet, voir le beau livre de Ne'eman et Kirsh\cite{Neeman}. Citons l'introduction de \cite{Fra} :{\it ``In
the 1960's, one of the mysteries in strong interaction physics was the
enormous proliferation of strongly interacting particles or
hadrons. Hadronic resonances seemed to exist with rather high
spin... the resonances were so numerous that it was not plausible that
they were all fundamental...''}.  Un exemple typique de telles
interactions est la d\'esint\'egration baryon-antibaryon (par exemple
neutron-antiproton) en trois m\'esons $\pi$ (\cite{Schwinger2}).

%\begin{center}
\begin{figure}[htbp]
\begin{center}
\begin{fmfgraph*}(100,100)
     \fmfpen{thick}
      \fmfleft{i1,i2}
      \fmfright{o1,o2,o3}
      \fmf{plain,label=$p$,label.side=left}{i1,v}
      \fmf{plain,label=$n$,label.side=left}{i2,v}
      \fmf{fermion}{i1,v}
      \fmf{fermion}{i2,v}
      \fmfpen{thin}
      \fmf{plain}{v,o1}
      \fmf{plain}{v,o2}
     	%\fmffreeze
      \fmf{plain}{v,o3}
      %\fmf{photon}{v3,v2}
   %   \fmffreeze
      %\fmf{photon,tension=0}{v3,o2}
%      \fmflabel{$\pi$}{i1}
   \end{fmfgraph*}
\caption{D\'esint\'egration proton-neutron}
\label{fig:mongraphe}
\end{center}
\end{figure}
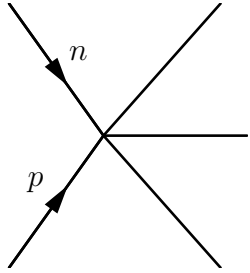

%   \end{center}

 L'analyse
complexe a permis de mettre au point des mod\`eles analytiques pour la th\'eorie; ces
mod\`eles pr\'esentent une sym\'etrie de jauge qui s'exprime en terme
d'op\'erateurs : les c\'el\`ebres op\'erateurs $L_n$ sont ainsi
apparus pour la premi\`ere fois en physique dans l'article de Virasoro
\cite{Vi1,Vi2} tr\`es souvent cit\'e mais rarement lu. Comme il est d'usage en physique, on ne voit pas appara\^itre le groupe ou l'alg\`ebre de Lie de fa\c con intrins\`eque, mais une repr\'esentation  donn\'ee en  g\'en\'eral par des op\'erateurs pour lesquels on doit ensuite v\'erifier les conditions de groupe ou d' alg\`ebre. Ici, ces op\'erateurs $L_n$
s'expriment comme somme quadratique d'oscillateurs, comme pour les
repr\'esentations  dans l'espace de Fock bosonique que nous avons
d\'efini plus haut. Dans \cite{Vi1} l'auteur conclut de fa\c con erron\'ee \`a l'existence d'une repr\'esentation lin\'eaire et non projective, donc le cocycle de Virasoro n'a pas \'et\'e d\'ecouvert par Miguel Virasoro lui-m\^ eme!  La propri\'et\'e d'alg\`ebre de Lie
sans le terme central, autrement dit $[L_n,L_m]=(m-n)L_{m+n}$, est quant \`a elle
apparue pour la premi\`ere fois dans un contexte physique dans
l'article de S. Fubini et G. Veneziano \cite{FuVe}. La mise en \'evidence du
terme central dans cette situation est attribu\'ee par P. Ramond \`a
Joe Weiss et Louis Clavelli (\cite{Ram}), voir aussi Brower et Thorn\cite{Brower}. Peu de temps apr\`es a \'et\'e construite
l'alg\`ebre de Neveu-Schwarz \cite{NeSc} avec la dualit\'e
bosons-fermions, qui devait marquer l'essor des superalg\`ebres...

Ce n'est pas notre but que de retracer toute cette
page glorieuse de l'histoire de la physique contemporaine, mais
mentionnons cependant que les premiers {\it mod\`eles de cordes et de
supercordes} sont issus de ces mod\`eles duaux: ce furent le mod\`ele
de Veneziano fond\'e sur les fonctions eul\'eriennes comme fonctions
d'amplitude (1968), puis la corde bosonique de Nambu (1970); la g\'eom\'etrie sous-jacente, avec l'action de ${\rm
Diff}(S^1)$ est apparue  explicitement pour la premi\`ere fois dans le travail de Galli\cite{Galli}; les
r\'esonances appara\^issant comme celles des vibrations de la
``corde''. Pour une approche \'el\'ementaire, mais aussi esth\'etique que p\'edagogique de la th\'eorie des cordes, voir\cite{Zwiebach}. La suite de l'histoire est longue, riche en rebondissements
passionnants, et n'est certainement pas achev\'ee de nos jours ; voir
encore les introductions de \cite{Zwiebach} et \cite{GreenSchwarzWitten}.

C'est en th\'eorie des cordes
qu'est apparue {\it l'anomalie conforme}, dite aussi
anomalie de Virasoro; pr\'ecisons maintenant pourquoi la g\'eom\'etrie conforme intervient dans ce sc\'enario. Tout le monde conna\^it la g\'eom\'etrie conforme et ses transformations qui pr\'eservent les angles, ainsi que le th\'eor\`eme de Liouville qui caract\'erise les transformations conformes en dimension $n\ge3$, et les hom\' eomorphismes biholomorphes si $n=2$. En th\'eorie des champs le groupe conforme intervient comme groupe de sym\'etries; l'approche axiomatique de la th\'eorie impose\cite{Wight} l'invariance sous l'action d'un groupe de sym\'etries de l'espace-temps, en g\'en\'eral celui de Poincar\'e; si on d\'ecide de l'\'etendre au groupe conforme, on obtient  la {\it th\'eorie des champs conformes}(CFT suivant le sigle anglais). La th\'eorie des cordes rend n\'ecessaire une th\'eorie quantique conforme sur l'espace-temps, d'o\`u cette anomalie de Virasoro:la th\'eorie des cordes bosoniques n'est invariante par
une transformation conforme de la m\'etrique que si la dimension de
l'espace-temps est \'egale \`a $26$. Le r\'esultat s'obtient par un
calcul explicite de la charge centrale (\cite{GreenSchwarzWitten}, p 130) pour
une certaine repr\'esentation de l'alg\`ebre de Virasoro; ce r\'esultat fut
montr\'e par Polyakov vers 1980 par des techniques d'int\'egrales de
chemin; ce fut la premi\`ere apparition du myst\'erieux $d=26$, dont
on peut rencontrer d'autres avatars dans diff\'erents contextes quantiques\cite{BowRaj}.

En mani\`ere de conclusion, nous allons expliquer la relation entre
l'alg\`ebre de Virasoro et \textit{l'effet Casimir}. Ce dernier est une
manifestation de l'\'energie du vide; il est classique que la
quantification canonique d'un oscillateur harmonique donne comme valeur du
hamiltonien 

\begin{equation}
\label{15}
\displaystyle H=(n+\frac{1}{2})\hslash \omega
\tag{15}
\end{equation} 

o\`u $\omega$ est la fr\'equence et $n$ le nombre de
particules, par cons\'equent $H \ne 0$ m\^eme si $n=0$. En 1948, H.B.G. Casimir
d\'eterminait la force d'attraction entre deux plaques conductrices parall\`eles,
par un calcul audacieux mais classique dans ce genre de physique :
l'\'energie de la ``bo\^\i te'' appara\^\i t comme la diff\'erence
entre deux quantit\'es infinies, l'\'energie du vide calcul\'ee
suivant la formule (12) et celle du vide perturb\'e par la
pr\'esence des plaques (voir \cite{Cas}). La confirmation exp\'erimentale
est due \`a Sparnaay en 1958 \cite{Sp}. Ce ph\'enom\`ene a trouv\'e une
autre interpr\'etation dans le cadre de la th\'eorie des champs
conformes en dimension $2$, dans le travail de Bl\"ote, Cardy et
Nightingale \cite{BlCaNi}, et c'est l\`a qu'intervient l'alg\`ebre de
Virasoro.

 Dans les th\'eories de champs conformes en dimension $2$,
les transformations du tenseur d'\'energie impulsion s'\'ecrivent sous
la forme

\begin{equation*}
\label{16}
\displaystyle T(z)dz^2 =T'(z')(dz')^2 + \frac{c}{12} \{z',z \} dz^2 
\tag{16}
\end{equation*}

o\`u $\{z',z\}$ d\'esigne la d\'eriv\'ee Schwarzienne de la
transformation $z'=\phi(z)$ (voir \cite{DrouffeItzykson1} Tome II p 103). Dans
le formalisme g\'eom\'etrique des champs de tenseurs, et de leurs transformations, cette relation s'\'ecrit
:

\begin{equation*}
\label{17}
\displaystyle \phi^* (Tdz^2)=[(T \circ \phi)(\phi')^2 + \frac{c}{12}
S(\phi) ]dz^2.
\tag{17}
\end{equation*}

Il ne s'agit que de la complexification de la formule de l'action
coadjointe de ${\rm Diff}(S^1)$, voir\cite[chap. 4]{Roger}.

%\mydessf{60}{60}{Cas.eps}{0}{Cas}{}

Pour l'\'etude de l'effet Casimir, on consid\`ere une bande de largeur
$L$ dans le plan complexe de la variable $u$ et la transformation
conforme $\displaystyle z=\exp \Big(\frac{2 i \pi u}{L}\Big)$ de cette
bande sur le plan des $z$. La formule (13) donne dans ce cas :
$$T_{({\rm bande})}(u)(du)^2=-\Big(\frac{2 \pi}{L}\Big)^2\Big(T(z)z^2
-\frac{c}{24}\Big)(dz)^2.$$
L'invariance par translation montre que les valeurs moyennes sur le
plan sont nulles donc $<T(z)>=0$. On en d\'eduit $\displaystyle
<T_{({\rm bande})}(u)>=\frac{c}{24}\Big(\frac{2\pi}{L}\Big)^2$ 
(\cite{DrouffeItzykson1}
p 105). Itzykson et Drouffe interpr\`etent ce r\'esultat de la fa\c
con suivante : {\it`` ... ceci nous permet d'interpr\'eter l'anomalie
comme un effet Casimir, c'est-\`a-dire un d\'eplacement de l'\'energie
libre comme cons\'equence de la g\'eom\'etrie finie...''}.

La comparaison avec la formule obtenue par Casimir \cite{Cas} est
instructive ; il serait int\'eressant de pouvoir faire directement le
lien entre l'\'energie de l'oscillateur et la valeur de la charge
centrale pour une repr\'esentation appropri\'ee de l'alg\`ebre de Virasoro.
%\begin{thebibliography}
%\bibdata{ARTHIST}
%\end{thebibliography}
\end{fmffile}
\bibliographystyle{plain}
\bibliography{ARTHIST}
\end{document}